\documentclass{amsart}
\usepackage{amssymb}
\usepackage{graphics}

\usepackage[
hyperref]{degt}

\let\PLUS+
{\catcode`\*\active \catcode`\[\active \catcode`\+\active
\let\+\relax
\gdef\deMaple{\catcode`\*\active \catcode`\[\active \catcode`\+\active
\def*{}\def[##1]{_{##1}}\def+{{}\PLUS{}}\catcode`\#\the\catcode`\%}}

\def\by{\bar y}
\def\nr{\subtext{nr}}
\def\trs{\subtext{tor}}

\let\slope\varkappa     %% slope of a curve
\def\tslope{\bar\slope} %% global slope
\let\lm\mu              %% local monodromy
\def\tlm{\tilde\lm}     %% twisted local monodromy
\let\pth\rho            %% paths from x_0 to x_i
\def\bm{\frak m}
\def\tbm{\tilde\bm}
\def\pim{\pth\subtext{im}}

\let\rra\Ga
\def\rrb{\Gb'}
\def\rrc{\Gb''}

\let\Mn=m  %% 'typical' Milnor number

\newcounter{line}

\def\lineref#1{\singset{#1}, \autoref{l.#1}}

\def\ttg{\text{\tt g}}
\def\tta{\text{\tt a}}
\def\nt#1{{\tt nt#1}}

\begin{document}

\title{On plane sextics with double singular points}

\author{Alex Degtyarev}
\address{Bilkent University\\
Department of Mathematics\\
06800 Ankara, Turkey}
\email{degt@fen.bilkent.edu.tr}



\begin{abstract}
We compute the fundamental groups of five maximizing sextics with double
singular points only; in four cases, the groups are as expected.
The approach used would apply to other sextics as well, given their
equations.
\end{abstract}

\maketitle

\section{Introduction}

The
fundamental group $\pi_1:=\pi_1(\Cp2\sminus D)$ of a plane curve
$D\subset\Cp2$, introduced by O.~Zariski in~\cite{Zariski:group},
is an important topological invariant of the curve.
Apart from distinguishing the connected components of the equisingular moduli
spaces, this group can be used as a seemingly inexpensive way of studying
algebraic surfaces, the curve serving as the branch locus of
a projection of the surface onto~$\Cp2$.

At present, the fundamental groups of all curves of degree up to five are
known, and the computation of the groups of irreducible curves of degree six
(sextics) is close to its completion, see~\cite{degt:book} for the principal
statements and further references.
In higher degrees, little is known: there are a few general theorems,
usually bounding the complexity of the group of a curve with sufficiently
`moderate' singularities, and a number of sporadic example scattered in the
literature. For further details on this fascinating subject, we refer the
reader to the recent
surveys~\cite{Artal:survey,Libgober:survey.2007,Libgober:survey2}.

\subsection{Principal results}
If a sextic $D\subset\Cp2$ has a singular point~$P$ of multiplicity three or
higher, then, projecting from this point, we obtain a trigonal (or, even
better, bi- or monogonal) curve in a Hirzebruch surface, see
\autoref{s.Hirzebruch}.
By means of
the so-called \emph{dessins d'enfants}, such curves and their topology can be
studied in purely combinatorial terms, as certain graphs in the plane.
The classification of such curves and the computation of their fundamental
groups were completed in~\cite{degt:book}.
If all singular points are double, the best that one can obtain is a
tetragonal curve, which is a much more complicated object.
(A reduction of tetragonal curves to trigonal curves in the presence of a
section is discussed in \autoref{s.resolvent}, see \autoref{rem.resolvent}.
It is the extra section that makes the problem difficult.)
At present, I do not know how the group of a tetragonal curve can be computed
unless the curve is real and its defining equation is known (and even then,
the approach
suggested in the paper
may still fail, \cf. \autoref{rem.A19}).

There is a special class of irreducible sextics, the so called
\emph{$\DG{2n}$-sextics} and, in particular, sextics of \emph{torus type}
(see \autoref{s.sextics} for the precise definitions),
for which the fundamental
group is non-abelian for some simple homological reasons,
see~\cite{degt:Oka}.
(The
fact that a sextic is of torus type is usually indicated by the presence of a
pair of parentheses in the notation; their precise meaning is explained in
\autoref{s.sextics}.)
On the other hand, thanks to the special structures and symmetries of these
curves, their explicit equations are known,
see~\cite{degt:Oka3,degt.Oka,Oka.Pho:moduli}.
In this paper, we almost complete the computation of the fundamental groups
of $\DG{2n}$-sextics (with one pair of complex conjugate sextics of torus
type left). Our principal results can be stated as follows.

\theorem\label{th.3A6+A1}
The fundamental group of the $\DG{14}$-special sextic with the set of
singularities \lineref{3A6+A1} in \autoref{tab.sextics} is
$\CG3\times\DG{14}$.
\endtheorem

\theorem\label{th.torus}
The fundamental groups of the irreducible sextics of torus type
with the sets of singularities
\lineref{(A14+A2)+A3},
\lineref{(A14+A2)+A2+A1},
and \lineref{(A11+2A2)+A4}
in \autoref{tab.sextics} are isomorphic to $\MG:=\CG2*\CG3$.
The group of the curve with the set of singularities
\lineref{(A8+3A2)+A4+A1} is
\begin{multline}
\pi_1=\bigl\<\Ga_2,\Ga_3,\Ga_4\bigm|
[\Ga_3,\Ga_4]=\{\Ga_2,\Ga_3\}_3=\{\Ga_2,\Ga_4\}_9=1,\\
\Ga_4\Ga_2\Ga_3\1\Ga_4\Ga_2\Ga_4(\Ga_4\Ga_2)^{-2}\Ga_3=
 (\Ga_2\Ga_4)^2\Ga_3\1\Ga_2\Ga_4\Ga_3\Ga_2\>,
\label{gr.(A8+3A2)+A4+A1}
\end{multline}
where $\{\Ga,\Gb\}_{2k+1}:=(\Ga\Gb)^k\Ga(\Ga\Gb)^{-k}\Gb\1$.
\endtheorem

\autoref{th.3A6+A1} is proved in \autoref{s.3A6+A1},
and \autoref{th.torus} is proved in
\autoref{s.(A14+A2)+A3}--\autoref{s.(A8+3A2)+A4+A1},
one curve at a time.
I do not know whether
the last group \eqref{gr.(A8+3A2)+A4+A1} is isomorphic to~$\MG$: all
`computable' invariants seem to coincide, see \autoref{rem.unknown},
but the presentations obtained resist all simplification attempt.
The quotient of~\eqref{gr.(A8+3A2)+A4+A1} by the extra relation
%$\Ga_3=\Ga_4$
$\{\Ga_2,\Ga_4\}_3=1$
is~$\MG$.

The next proposition is proved in \autoref{proof.pert}.
(The perturbation $\singset{3A6+A1}\to\singset{3A6}$ excluded in the
statement results in a $\DG{14}$-special sextic and the fundamental group
equals $\CG3\times\DG{14}$, see~\cite{degt.Oka}.)

\proposition\label{prop.pert}
Let~$D'$ be a nontrivial perturbation of a sextic as in
Theorems~\ref{th.3A6+A1} or~\ref{th.torus}.
Unless the set of singularities of~$D'$ is \singset{3A6},
the group $\pi_1(\Cp2\sminus D')$ is $\MG$ or~$\CG6$, depending on
whether $D'$ is or, respectively, is not of torus type.
\endproposition

With \autoref{th.3A6+A1} in mind,
the fundamental groups of all $\DG{2n}$-special sextics, $n\ge5$, are known,
see~\cite{degt:book}.
Modulo the feasible conjecture that any sextic of torus type degenerates to a
maximizing one,
the only such sextic whose group remains unknown is
\lineref{(A8+A5+A2)+A4} in \autoref{tab.sextics}.
(This
conjecture has been proved, and \emph{all} groups except the one just
mentioned are indeed known; details will appear elsewhere.)
Most of these groups are isomorphic to~$\MG$, see~\cite{degt:book} for
details and further references.

I would like to mention
an
alternative approach, see~\cite{Artal:trends}, reducing a plane sextic with
large Milnor number to a trigonal curve equipped with a
number of sections, all but one splitting in the covering elliptic surface.
It was used in~\cite{Artal:trends} to handle the curves in
lines~\ref{l.A19}--\ref{l.A15+A4} in \autoref{tab.sextics}.
This approach is also used in a forthcoming paper to produce the defining
equations of most sextics listed in \autoref{tab.sextics}; then, the
fundamental groups of most \emph{real} ones can be computed using
\autoref{th.vanKampen}. All groups that could be found are abelian.
Together with the classification of sextics, which is also almost completed,
this fact implies that, with very few exceptions, the fundamental group of a
non-special irreducible simple sextic is abelian.

\subsection{Idea of the proof\noaux{ (see \autoref{s.strategy} for more details)}}\label{s.idea}
We use
the classical Zariski--van Kampen method, \cf. \autoref{th.vanKampen},
expressing the fundamental group of a curve in terms of its braid monodromy
with respect to an appropriate pencil of lines.
The curves and pencils considered are \emph{real}, and the braid monodromy in
a neighborhood of the real part of the pencil is computed in terms of the
real part of the curve. (This approach originates in topology of real
algebraic curves; historically,
it goes back to Viro, Fiedler, Kharlamov, Rokhlin, and Klein.)
Our main contribution is the description of the monodromy
along a real segment where all four branches of the curve are non-real, see
\autoref{prop.gamma}. Besides, the curves are not required to be
\emph{strongly real}, \ie, non-real singular fibers are allowed.
Hence, we follow Orevkov \cite{Orevkov:braids} and
attempt to extract information
about such non-real fibers from the real part of the curve.
The outcome is \autoref{th.vanKampen},
which gives us
an `upper bound' on the
fundamental group in question.
The applicability issues and a few other common tricks are
discussed in \autoref{s.strategy}.

\subsection{Contents of the paper}
In \autoref{S.sextics}, we introduce the terminology related to plane
sextics, list the sextics that are still to be investigated, and discuss
briefly the few known results.
In \autoref{S.monodromy}, we outline an approach to the (partial) computation
of the braid monodromy of a real tetragonal curve and state an appropriate
version of the Zariski--van Kampen theorem.
Finally, in \autoref{S.proofs} the results of \autoref{S.monodromy}
and known equations are used
to prove Theorems~\ref{th.3A6+A1} and~\ref{th.torus} and \autoref{prop.pert}.

\subsection{Conventions}
All group actions are \emph{right}. Given a right action
$X\times G\to X$ and a pair of elements $x\in X$, $g\in G$, the image of
$(x,g)$ is denoted by $x\ra g\in X$.
The same postfix notation and multiplication convention is
often
used for maps: it is under this convention that the monodromy
$\pi_1(\text{base})\to\Aut(\text{fiber})$
of a locally trivial fibration is a
homomorphism rather than an anti-homomorphism.

The assignment symbol $:=$ is used as a shortcut for `is defined as'.

We use the conventional symbol \qedsymbol\ to mark the ends of the proofs.
Some statements are marked with \donesymbol\ or \pnisymbol: the former means
that the proof has already been explained (for example, most corollaries),
and the latter indicates that the proof is not found in the paper and the
reader is directed to the literature, usually cited at the beginning of the
statement.

\subsection{Acknowledgements}
This paper was written during my sabbatical stay at
\emph{l'Instutut des Hautes \'{E}tudes Scientifiques} and
\emph{Max-Planck-Institut f\"{u}r Mathematik};
I would like to thank these institutions for their support and hospitality.
I am also grateful to M.~Oka, who kindly clarified for me the results
of~\cite{Oka.Pho:moduli},
to V.~Kharlamov, who patiently introduced me to the history of the subject,
and to the
anonymous
referee of the paper, who made
a number of valuable suggestions
and checked
and confirmed the somewhat unexpected result of \autoref{s.(A8+3A2)+A4+A1}.

\section{Preliminaries}\label{S.sextics}

\subsection{Special classes of sextics}\label{s.sextics}
A plane sextic $D\in\Cp2$ is called \emph{simple} if all its singularities
are simple, \ie, those of type $\bA$--$\bD$--$\bE$.
The total Milnor number~$\mu$ of a simple sextic~$D$ does not exceed~$19$,
see~\cite{Persson:sextics}; if $\mu=19$, then $D$ is called
\emph{maximizing}.
Maximizing sextics are always defined over algebraic number fields and
their moduli spaces are discrete: two such sextics are
equisingular deformation equivalent if and only if they are related by a
projective transformation of~$\Cp2$.

A sextic~$D$ is said to be of \emph{torus type} if its equation can be
represented in the form $f_2^3+f_3^2=0$, where $f_2$ and~$f_3$ are some
polynomials of degree~$2$ and~$3$, respectively.
The points of intersection of the conic $\{f_2=0\}$ and cubic~$\{f_3=0\}$ are
always singular for~$D$. These singular points play a very special r\^{o}le;
they are called the \emph{inner} singularities (with respect to the given
torus structure).
For the vast majority of curves, a torus
structure is unique, and in this case it is
common to parenthesize the inner singularities in the notation.

An irreducible sextic~$D$ is called \emph{$\DG{2n}$-special} if its
fundamental group $\pi_1(\Cp2\sminus D)$ admits a dihedral quotient
$\DG{2n}:=\CG{n}\rtimes\CG2$.
According to~\cite{degt:Oka}, only $\DG6$-, $\DG{10}$-, and $\DG{14}$-special
sextics exist, and an irreducible sextic is of torus type if and only if it
is $\DG6$-special. (In particular, torus type is a topological property.)

Any sextic~$D$
of torus type is a degeneration of Zariski's six-cuspidal sextic, which is
obtained from a generic pair $(f_2,f_3)$.
It follows that the fundamental group of~$D$
factors to the modular group
$\Gamma:=\SL(2,\Z)=\CG2*\CG3=\BG3/(\Gs_1\Gs_2\Gs_1)^2$,
see~\cite{Zariski:group}; in particular, this group is infinite.
Conjecturally, the fundamental groups of all other irreducible simple sextics are
finite.

\subsection{Sextics to be considered}\label{s.table}
It is expected that, with few explicit exceptions (\eg, \singset{9A2}),
any simple sextic
degenerates to a maximizing one. (The proof of this conjecture,
which relies upon the theory of $K3$-surfaces, is currently
a work in progress. In fact, most curves degenerate to one of those whose
group is already known.)
Hence, it is essential to compute the fundamental groups of the maximizing
sextics; the others would follow.
The groups of all irreducible sextics with a singular point of multiplicity
three or higher are known, see~\cite{degt:book} for a summary of the results,
and those with $\bA$ type singularities only are still to be investigated.

\table
\caption{Irreducible maximizing sextics with $\bA$ type singularities}\label{tab.sextics}
\def\torus{(torus type)}
\def\spec#1{($\DG{#1}$-sextic)}
\def\*{\rlap{$^*$}}
\def\1{(1,0)}
\def\countr{0}
\def\countc{0}
\def\savecount#1#2{\count0=#1 \advance\count0 #2\xdef#1{\the\count0}}
\def\savecounts(#1,#2){\savecount\countr#1\savecount\countc#2$(#1,#2)$}
\hbox to\hsize{\hss\vbox{\obeylines\let
\cr%
\halign{\gdef\LABEL{#}\space\strut\hss%
 \refstepcounter{line}\theline.\label{l.\LABEL}%
 &\quad\singset{\LABEL}#\hss%
% &\quad\hss$#$\hss%
 &\quad\hss\expandafter\savecounts#\hss%
 &\quad#\hss\space
\noalign{\hrule\vspace{3pt}}%
\omit\strut\hss$\#$\ &\omit\quad Singularities\hss&\omit\quad\hss$(r,c)$\hss%
 &Equation, $\pi_1$, remarks
\noalign{\vspace{2pt}\hrule\vspace{2pt}}%
A19&&(2,0)&$\pi_1=\CG6$, see \cite{Artal:trends}
A18+A1&&(1,1)&$\pi_1=\CG6$, see \cite{Artal:trends}
(A17+A2)&&(1,0)\*&$\pi_1=\MG$, see \cite{Artal:trends,degt:2a8} \torus
A16+A3&&(2,0)&$\pi_1=\CG6$, see \cite{Artal:trends}
A16+A2+A1&&(1,1)&$\pi_1=\CG6$, see \cite{Artal:trends}
A15+A4&&(0,1)\*&$\pi_1=\CG6$, see \cite{Artal:trends}
A14+A4+A1&&(0,3)&
(A14+A2)+A3&&\1&$\pi_1=\MG$, see \autoref{s.(A14+A2)+A3} \torus
(A14+A2)+A2+A1&&\1&$\pi_1=\MG$, see \autoref{s.(A14+A2)+A2+A1} \torus
A13+A6&&(0,2)
A13+A4+A2&&(2,0)
A12+A7&&(0,1)
A12+A6+A1&&(1,1)
A12+A4+A3&&\1
A12+A4+A2+A1&&(1,1)
A11+2A4&&(2,0)
(A11+2A2)+A4&&\1&$\pi_1=\MG$, see \autoref{s.(A11+2A2)+A4} \torus
A10+A9&&(2,0)\*
A10+A8+A1&&(1,1)
A10+A7+A2&&(2,0)
A10+A6+A3&&(0,1)
A10+A6+A2+A1&&(1,1)
A10+A5+A4&&(2,0)
A10+2A4+A1&&(1,1)
A10+A4+A3+A2&&\1
A10+A4+2A2+A1&&(2,0)
A9+A6+A4&&(1,1)\*&
A9+2A4+A2&&(1,0)\*&$\pi_1=\eqref{eq.A9+2A4+A2}$, see \cite{degt:Oka3} \spec{10}
(2A8)+A3&&\1&$\pi_1=\MG$, see \cite{degt:2a8} \torus
A8+A7+A4&&(0,1)
A8+A6+A4+A1&&(1,1)
(A8+A5+A2)+A4&&(0,1)&\nt{104} in~\cite{Oka.Pho:moduli} \torus
(A8+3A2)+A4+A1&&\1&$\pi_1=\eqref{gr.(A8+3A2)+A4+A1}$, see \autoref{s.(A8+3A2)+A4+A1} \torus
A7+2A6&&(0,1)
A7+A6+A4+A2&&(2,0)
A7+2A4+2A2&&\1
3A6+A1&&\1&$\pi_1=\CG3\times\DG{14}$, see \autoref{s.3A6+A1} \spec{14}
2A6+A4+A2+A1&&(2,0)
A6+A5+2A4&&(2,0)
\noalign{\vspace{2pt}\hrule\vspace{3pt}\footnotesize%
\hbox{\strut\quad%
Marked with a $^*$ are the sets of singularities realized by reducible %
sextics as well}%
\hbox{\quad%
There are $\countr$ real and $\countc$ pairs of complex conjugate curves}%
}\crcr}}\hss}
\endtable

A list of irreducible maximizing sextics with $\bA$ type singular points only
can be compiled using the results of~\cite{Yang} (a list of the sets of
singularities realized by such sextics)
and~\cite{Shimada:maximal} (a description of the moduli spaces).
We represent the result in \autoref{tab.sextics}, where the column $(r,c)$
shows the number of classes: $r$ is the number of real sextics, and $c$ is
the number of pairs of complex conjugate ones.
The approach developed further in the paper lets one compute (or
at least estimate)
the fundamental group of a sextic with $\bA$ type singularities, provided that
its equation is known. In the literature, I
could find explicit equations for lines~\ref{l.A19}--\ref{l.A15+A4},
\ref{l.(A14+A2)+A3}, \ref{l.(A14+A2)+A2+A1}, \ref{l.(A11+2A2)+A4},
\ref{l.A9+2A4+A2}, \ref{l.(2A8)+A3}, \ref{l.(A8+A5+A2)+A4},
\ref{l.(A8+3A2)+A4+A1}, and~\ref{l.3A6+A1}.
With the results of this paper
(Theorems~\ref{th.3A6+A1} and~\ref{th.torus}) taken into account, the
groups of all these sextics except \lineref{(A8+A5+A2)+A4} (which is not
real) are known.

\remark\label{rem.A19}
Unfortunately, our approach does not always work even if the curve is real.
Thus, each of the two sextics with the set of singularities \lineref{A19}
has a single real point (the isolated singular point of type~$\bA_{19}$;
see~\cite{Artal:trends} for the equations) and \autoref{th.vanKampen} does
not provide enough relations to compute the group.
\endremark

\subsection{Known results}
The fundamental group of the $\DG{10}$-special sextic with the set of
singularities \lineref{A9+2A4+A2} in \autoref{tab.sextics},
can be described as follows, see \cite{degt:Oka3}
(where $'$ temporarily stands for the commutant of a group):
\[
\pi_1/\pi_1''=\CG3\times\DG{10},\quad
\pi_1''=\SL(2,\Bbbk_9),
\label{eq.A9+2A4+A2}
\]
where $\Bbbk_9$ is the field of nine elements.
The fundamental groups of the first twelve sextics,
lines~\ref{l.A19}--\ref{l.A15+A4}, have been found in~\cite{Artal:trends}:
with the exception of \lineref{(A17+A2)} (sextic of torus type, $\pi_1=\MG$),
they are all abelian.
To my knowledge, the groups not mentioned in \autoref{tab.sextics} have not
been computed yet.

\section{The braid monodromy}\label{S.monodromy}

\subsection{Hirzebruch surfaces}\label{s.Hirzebruch}
A \emph{Hirzebruch surface} $\Sigma_d$, $d>0$, is a geometrically ruled
rational surface with a
(unique) \emph{exceptional section}~$E$ of
self-intersection~$-d$.
Typically, we use affine coordinates $(x,y)$ in~$\Sigma_d$ such that $E$ is
given by $y=\infty$; then, $x$ can be regarded as an affine coordinate in the
base of the ruling.
(The
line $\{x=\infty\}$ plays no special r\^{o}le; usually, it is assumed
sufficiently generic.)
The fiber of the ruling over a point~$x$ in the base is denoted by~$F_x$, and
the \emph{affine fiber} over~$x$ is $F_x^\circ:=F_x\sminus E$.
This is an affine space over~$\C$; in particular, one can speak about convex
hulls in $F_x^\circ$.

An \emph{$n$-gonal curve} is a reduced curve $C\subset\Sigma_d$ intersecting
each fiber at $n$ points, \ie, such that the restriction to~$C$ of the ruling
$\Sigma_d\to\Cp1$ is a map of degree~$n$. A \emph{singular fiber} of an
$n$-gonal curve~$C$ is a fiber~$F$ of the ruling intersecting $C+E$
geometrically
at fewer than $(n+1)$ points.
A singular fiber~$F$ is \emph{proper} if $C$ does not pass through $F\cap E$.
The curve~$C$ is \emph{proper} if so are all its singular fibers.
In
other words, $C$ is proper if it is disjoint from~$E$.

In affine coordinates $(x,y)$ as above an $n$-gonal curve $C\subset\Sigma_d$
is given by a polynomial of the form $\sum_{i=0}^na_i(x)y^i$, where
$\deg a_i\le m+d(n-i)$ for some $m\ge0$ (in fact, $m=C\cdot E$) and at least
one polynomial~$a_i$ does have the prescribed degree (so that $C$ does not
contain the fiber $\{x=\infty\}$).
The curve is proper if and only if $m=0$; in this case $a_n(x)=\const$.

A proper $n$-gonal curve $C\subset\Sigma_d$ defines a distinguished
\emph{zero section} $Z\subset\Sigma_d$, sending each point $x\in\Cp1$ to the
barycenter of the $n$ points of $F_x^\circ\cap C$.
Certainly,
this section does not need to coincide with $\{y=0\}$, which
depends on the choice of the coordinates.

\subsection{The cubic resolvent}\label{s.resolvent}
Consider a reduced \emph{real} quartic polynomial
\[
f(x,y):=y^4+p(x)y^2+q(x)y+r(x),
\label{eq.f}
\]
so that its roots
$y_1,y_2,y_3,y_4$ (at each point~$x$)
satisfy $y_1+y_2+y_3+y_4=0$, and consider
the \emph{\rom(modified\rom) cubic resolvent} of~$f$
\[
y^3-2p(x)y^2+b_1(x)y+q(x)^2,\qquad
b_1:=p^2-4r,
\label{eq.resolvent}
\]
and its reduced form
\[
\by^3+g_2(x)\by+g_3(x)
\label{eq.reduced}
\]
obtained by the substitution $y=\by+\frac23p$.
The discriminants of~\eqref{eq.f}--\eqref{eq.reduced}
are equal:
\[
\deMaple
D=16*p^4*r-4*p^3*q^2-128*p^2*r^2+144*p*q^2*r-27*q^4+256*r^3.
\label{eq.discrim}
\]
Recall that $D=0$ if and only if \eqref{eq.f} or, equivalently,
\eqref{eq.resolvent} or~\eqref{eq.reduced}
has a multiple root. Otherwise,
$D<0$ if and only if exactly two roots of~\eqref{eq.f} are real.
The roots of~\eqref{eq.resolvent} are
\[
\aligned
\Ga&:=(y_1+y_2)(y_3+y_4)=-(y_1+y_2)^2,\\
\Gb&:=(y_1+y_3)(y_2+y_4)=-(y_1+y_3)^2,\\
\Gg&:=(y_1+y_4)(y_2+y_3)=-(y_1+y_4)^2,
\endaligned
\label{eq.roots}
\]
and those of~\eqref{eq.reduced} are obtained from~\eqref{eq.roots}
by shifting the barycenter $\frac13(\Ga+\Gb+\Gg)$ to zero.

\remark\label{rem.resolvent}
If $\{f(x,y)=0\}$ is a proper tetragonal curve in a Hirzebruch
surface $\Sigma_d$, then \eqref{eq.resolvent} defines a proper trigonal curve
$C'\subset\Sigma_{2d}$
and a distinguished section $S:=\{y=0\}$
(in general, \emph{other} than the zero section)
which is tangent (more precisely, has even intersection
index at each intersection point) to~$C'$.
Conversely, \eqref{eq.f} can be recovered from~\eqref{eq.resolvent} (together
with the section $S=\{y=0\}$) uniquely up to the automorphism $y\mapsto-y$, which
takes~$q$ to $-q$.
\endremark

\remark\label{rem.q}
One has
\[*
q=-(y_1+y_2)(y_1+y_3)(y_1+y_4);
\]
hence, $q$ vanishes if and only if two of the roots of~\eqref{eq.f}
are opposite.
If all roots are non-real,
$y_{1,2}=\rra\pm\rrb i$, $y_{3,4}=-\rra\pm\rrc i$,
$\rra,\rrb,\rrc\in\R$, then
\[*
%q=2\rra(\rrb-\rrc)(\rrb+\rrc),\quad
q=2\rra\bigl({\rrb}^2-{\rrc}^2\bigr),\quad
b_1=-8\rra^2\bigl({\rrb}^2+{\rrc}^2\bigr)+\bigl({\rrb}^2-{\rrc}^2\bigr)^2.
%\label{eq.q,b1}
\]
Hence, $q(x)=0$ if and only if either $\rra=0$
(and then $b_1(x)>0$, assuming that $D(x)\ne0$)
or $\rrb=\pm\rrc$ (and then $b_1(x)<0$).
If $y_1=y_2$, \ie, $\rrb=0$, then $q(x)>0$ if and only if one has the inequality
$y_1<\Re y_3=\Re y_4$
equivalent to $y_1<0$.
\endremark

\remark\label{rem.g3}
Observe also that, if $y_1=y_2$, then $g_3$ takes the form
\[*
g_3=\frac2{27}(y_1-y_4)^3(y_1-y_3)^3.
%\label{eq.g3}
\]
Hence, $g_3(x)<0$ if and only if the two other roots are real and separated by
the double root $y_1=y_2$.
Otherwise, either $y_1<\Re y_3,\Re y_4$ or $y_1>\Re y_3,\Re y_4$, and, in
view of \autoref{rem.q}, the former holds if and only if $q(x)>0$.
\endremark

\subsection{The real monodromy}\label{s.settings}
Choose affine coordinates $(x,y)$ in the Hirzebruch surface~$\Sigma_d$ so
that the exceptional section~$E$ is $\{y=\infty\}$.
Consider a \emph{real} proper tetragonal curve $C\subset\Sigma_d$; it is
given by a real polynomial $f(x,y)$ as in~\eqref{eq.f}.
Over a generic real point $x\in\R$, the four points $y_1,\ldots,y_4$
of
the intersection
$C\cap F_x^\circ$ can be ordered lexicographically,
according to the \emph{decreasing} of $\Re y$ first and $\Im y$ second.
We always assume this ordering. Then, choosing a real reference point $y\gg0$,
we have a \emph{canonical} geometric basis $\{\Ga_1,\ldots,\Ga_4\}$
for the fundamental group
$\pi(x):=\pi_1(F_x^\circ\sminus C, y)$, see \autoref{fig.basis}.
\figure[ht]
\centerline{\cpic{basis}}
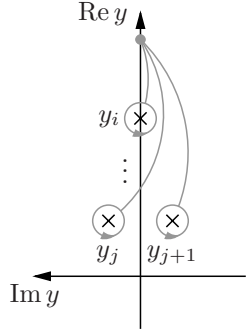
\caption{The canonical basis}\label{fig.basis}
\endfigure

Let $x_1,\ldots,x_r$ be all \emph{real} singular fibers of~$C$, ordered by
increasing. For each~$i$, consider a pair of nonsingular fibers
$x_i^-:=x_i-\epsilon$ and $x_i^+:=x_i+\epsilon$, where $\epsilon$ is a
sufficiently small positive real number, see \autoref{fig.base}.
\figure[ht]
\centerline{\cpic{base}}
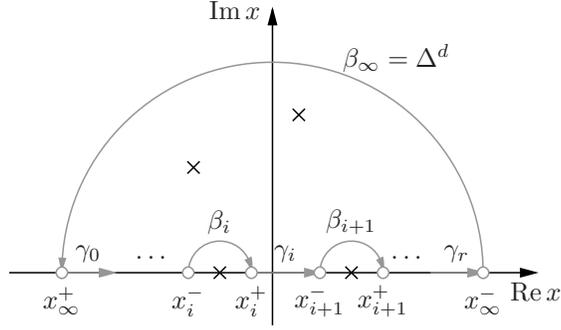
\caption{The monodromies $\Gb_i$ and $\Gg_j$}\label{fig.base}
\endfigure
Denote $x_0=x_{r+1}=\infty$ and,
assuming the fiber $x=\infty$ nonsingular, pick also a pair of
real
nonsingular fibers $x_{r+1}^-=x_\infty^-:=R\gg0$ and
$x_0^+=x_\infty^+:=-R$.
Identify all groups $\pi(x_i^\pm)$ with
the free group~$\FG4$ by means of their respective
canonical bases. (All reference points are chosen in a real section
$y=\const\gg0$, which is assumed disjoint from
the fiberwise convex hull of~$C$ over the disk
$\ls|x|\le R$.)
Consider the semicircles
$t\mapsto x_i+\epsilon e^{i\pi(1-t)}$, $t\in[0,1]$, and the line
segments $t\mapsto t$, $t\in[x_j^+,x_{j+1}^-]$, \cf.
\autoref{fig.base}.
These paths give rise to the \emph{monodromy isomorphisms}
\[*
\Gb_i\:\pi(x_i^-)\to\pi(x_i^+),\quad
\Gg_j\:\pi(x_j^+)\to\pi(x_{j+1}^-),
%\label{eq.monodromy}
\]
$i=1,\ldots,r$, $j=0,\ldots,r$.
In addition, we also have the monodromy
$\Gb_0=\Gb_\infty=\Gb_{r+1}\:\pi(x_\infty^-)\to\pi(x_\infty^+)$ along the
semicircle $t\mapsto Re^{i\pi t}$, $t\in[0,1]$,
and the \emph{local monodromies}
\[*
\lm_i\:\pi(x_i^+)\to\pi(x_i^+),\quad
i=1,\ldots,r
%\label{local}
\]
along the circles $t\mapsto x_i+\epsilon e^{2\pi it}$, $t\in[0,1]$.
Using the identifications $\pi(x_i^\pm)=\FG4$ fixed above, all
$\Gb_i$, $\lm_i$, $\Gg_j$
can be regarded as elements of the automorphism group $\Aut\FG4$, and
as such they belong to the braid group~$\BG4$.
Recall, see~\cite{Artin}, that Artin's \emph{braid group}
$\BG4\subset\Aut\<\Ga_1,\ldots,\Ga_4\>$ is the subgroup consisting of the
automorphisms taking each generator~$\Ga_i$ to a conjugate of a generator and
preserving the product $\Ga_1\Ga_2\Ga_3\Ga_4$.
It is generated by the three braids
\[*
\Gs_i\:\
\Ga_i\mapsto\Ga_i\Ga_{i+1}\Ga_i\1,\quad
\Ga_{i+1}\mapsto\Ga_i,\quad
i=1,2,3,
\]
the defining relations being
$\{\Gs_1,\Gs_2\}_3=\{\Gs_2,\Gs_3\}_3=[\Gs_1,\Gs_3]=1$.

\subsection{The computation}\label{s.monodromy}
The braids $\Gb_i$, $\lm_i$, and~$\Gg_j$ introduced in the previous
section are easily computed from the real part $C_\R\subset\R^2$ of the
curve. In the figures, we use the following notation:
\roster*
\item
real branches of~$C$ are represented by solid bold lines;
\item
pairs $y_i,y_{i+1}$ of complex conjugate branches are represented by dotted
lines (showing the common real part $\Re y_i=\Re y_{i+1}$);
\item
relevant fibers of~$\Sigma_d$ are represented by vertical dotted grey lines.
\endroster
Certainly, the dotted lines are not readily seen in the figures; however, in
most cases, it is only the intersection indices that matter, and the latter
are determined by the indexing of the branches at the starting and ending
positions.

We summarize the results in the next three statements.
The first one is obvious: essentially, one speaks about the link
of the singularity $y^4-x^{4d}$.

%\lemma[see, \eg, \cite{degt:book}]\label{lem.infty}
\lemma\label{lem.infty}
Assume that $R\gg0$ is so large that the disk $\{\ls|x|<R\}$ contains all
singular fibers of~$C$.
Then one has $\Gb_\infty=\Delta^d$, where
$\Delta:=\Gs_1\Gs_2\Gs_3\Gs_1\Gs_2\Gs_1\in\BG4$ is the
\emph{Garside element}.
\done
\endlemma

The following lemma is easily proved by considering the local normal forms of
the singularities.
(In
the simplest case of a vertical tangent, the circumventing braids~$\Gb$
are computed, \eg, in~\cite{Orevkov:braids}; the general case is completely
similar.)
For the statement, we extend the standard notation
$\bA_\Mn$, $\Mn\ge1$,
to $\bA_0$ to designate a simple tangency of~$C$ and the fiber.

\lemma\label{lem.beta}
The braids $\Gb_j$ and $\lm_j$ about a singular fiber~$x_j$ of
type~$\bA_\Mn$, $\Mn\ge0$, depend only on~$\Mn$ and the pair $(i,i+1)$ of
indices of the branches that merge at the singular point. They are as shown
in \autoref{fig.beta}.
\figure[ht]
\centerline{\cpic{A2k-1}\hfil\cpic{A2k-}\hfil\cpic{A2k+}}
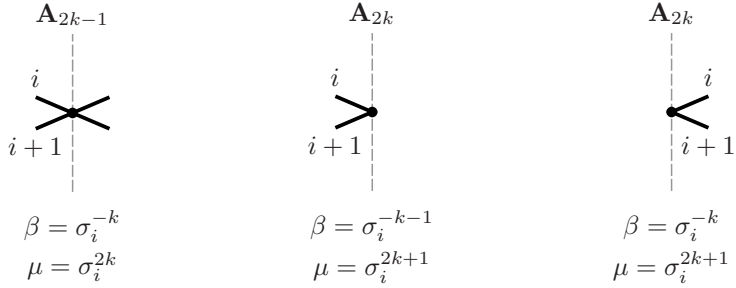
\caption{The braids $\Gb$ and $\lm$}\label{fig.beta}
\endfigure
\done
\endlemma

\remark
At a point of type $\bA_{2k-1}$, it is
\emph{not} important whether the two branches
of~$C$ at this point are real or complex conjugate. On the other hand, at a
point of type~$\bA_{2k}$ it \emph{does} matter whether the number of real
branches increases or decreases.
If a fiber contains two double points, with indices $(1,2)$ and $(3,4)$, then
the powers of~$\Gs_1$ and~$\Gs_3$ contributed to~$\Gb$ or $\lm$ by each
of the points are multiplied;
since $\Gs_1$ and $\Gs_3$ commute, the order is not important.
\endremark

The following statement is our principal technical tool,
most important being \autoref{fig.gamma}, right,
%which describes
describing
the behaviour of the `invisible' branches.
(Note that the two dotted lines in the figure may cross; the
permutation of the branches depends on the parity of the twist parameter~$t$
introduced in the statement.)

\proposition\label{prop.gamma}
Let $I$ be a real segment in the $x$-axis free of singular fibers of~$C$.
Then the monodromy~$\Gg$ over~$I$ is
\roster*
\item
identity, if all four branches of~$C$ over~$I$ are real, and
\item
as shown in \autoref{fig.gamma} otherwise.
\endroster
\figure[ht]
\centerline{\cpic{real-}\hfil\cpic{real+}\hfil\cpic{complex}}
\caption{The braids $\Gg$}\label{fig.gamma}
\endfigure
Here,
$\tau:=\Gs_2\1\Gs_3\Gs_1\1\Gs_2$ and
the \emph{twist parameter} $t$ in \autoref{fig.gamma}, right is the
number of roots $x'\in I$ of the coefficient~$q(x)$,
see \eqref{eq.f},
such that $b_1(x')>0$,
%\rom(see \eqref{eq.f}, \eqref{eq.resolvent}\rom)
see \eqref{eq.resolvent},
and $q$ changes sign at~$x'$\rom;
each root~$x'$ contributes~$+1$ or~$-1$ depending on
whether $q$ is increasing or decreasing at~$x'$, respectively.
\endproposition

\proof
The only case that needs consideration, \viz. that of four non-real
branches, see \autoref{fig.gamma}, right, is given by \autoref{rem.q}.
%(Note that the two dotted lines in the figure may cross; the resulting
%permutation depends on the parity of~$t$.)
Indeed, the canonical basis in the fiber $F_x^\circ$ over $x\in I$ changes
when the real parts of all four branches vanish, and this happens when
$q(x)=0$ and $b_1(x)>0$. This change contributes
$\tau^{\pm1}$ to~$\Gg$, and the sign $\pm1$ (the
direction of rotation) depends on whether $q$ increases or
decreases.
\endproof

\remark\label{rem.gamma}
A longer segment~$I$ with exactly two real branches of~$C$ over it can be
divided into smaller pieces $I_1,I_2,\ldots$, each containing a single
crossing point as in \autoref{fig.gamma}; then, the monodromy~$\Gg$ over~$I$
is the product of the contributions of each piece.
In fact, as explained above, the precise position and number of crossings is
irrelevant; what only matters is the final permutation between the endpoints
of~$I$. For example, to minimize the number of elementary pieces, one can
always assume the branches, both bold and dotted, monotonous.
\endremark

\subsection{The Zariski--van Kampen theorem}
We are interested in the fundamental group
$\pi_1:=\pi_1(\Sigma_{\tilde d}\sminus(\tilde C\cup E))$, where
$\tilde C\subset\Sigma_{\tilde d}$ is a real tetragonal curve, possibly
improper, and $E\subset\Sigma_{\tilde d}$ is the exceptional section.
To compute~$\pi_1$, we consider the \emph{proper model}
$C\subset\Sigma_d$, obtained from~$\tilde C$ by blowing up all points of
intersection $\tilde C\cap E$ and blowing down the corresponding fibers.
In addition to the braids~$\Gb_i$, $\lm_i$, and~$\Gg_j$ introduced in
\autoref{s.settings}, to each real singular fiber~$x_i$ of~$C$ we assign its
\emph{local slope} $\slope_i\in\pi(x_i^+)$,
which depends on the type of the
corresponding singular fiber of the original curve~$\tilde C$.
Roughly, consider a small analytic disk $\Phi\subset\Sigma_d$ transversal to
the fiber~$F_{x_i}$ and disjoint from~$C$ and~$E$, and a similar disk
$\tilde\Phi\subset\Sigma_{\tilde d}$ with respect to $\tilde C$.
Let $\tilde\Phi'\subset\Sigma_d$ be the image of~$\tilde\Phi$, and assume
that the boundaries $\partial\Phi$ and $\partial\tilde\Phi'$ have a common
point in the fiber over~$x_i^+$.
Then the loop $[\partial\tilde\Phi']\cdot[\partial\Phi]\1$ is homotopic to a
certain class $\slope_i\in\pi(x_i^+)$, well defined up to a few moves
irrelevant in the sequel. This class is the slope.

Roughly, the slope measures
(in the form of the twisted monodromy, see the definitions prior to
\autoref{th.vanKampen})
the deviation
of the
braid monodromy of an improper curve $\tilde C$ from that of its proper
model~$C$. Slopes appear in the relation at infinity as well, compensating
for the fact that, near improper singular fibers, the curve intersects
\emph{any} section of~$\Sigma_{\tilde d}$.
Details and further
properties are found in \cite[\S5.1.3]{degt:book}; in this paper, slopes
are used in \autoref{th.vanKampen}.

\remark\label{rem.slopes}
In all examples considered below,
$\tilde C\subset\Sigma_{d-1}$ has a single improper fiber~$F$,
where $\tilde C$ has a singular point of type~$\tA_\Mn$, $\Mn\ge1$, maximally
transversal to both~$E$ and~$F$.
If $F=\{x=0\}$, such a curve~$\tilde C$ is given by a
polynomial~$\tilde f$ of the form $\sum_{i=0}^4y^ia_i(x)$ with $a_4(x)=x^2$
and
%$a_3(0)=0$,
$x\mathrel|a_3(x)$,
and the defining polynomial of its transform $C\subset\Sigma_d$ is
$f\nr(x,y):=x^2\tilde f(x,y/x)$.
The corresponding singular fiber of~$C$ has a node $\bA_1$ at $(0,0)$ and
another double point $\bA_{\Mn-2}$ (assuming $\Mn\ge2$).
\endremark

Thus,
the only nontrivial example relevant in the sequel is the one
%given by \autoref{ex.slope}
described
below. (By the very definition,
%By definition,
at each singular fiber~$x_i$ proper for~$\tilde C$
the slope
is $\slope_i=1$.)
A great deal of other
examples of both computing the slopes
and using them in the study of the fundamental group
are found in~\cite{degt:book}.

\example\label{ex.slope}
At
the only improper fiber $x_i=0$ described in \autoref{rem.slopes}
the
slope is the class of $\Ga_j\Ga_{j+1}$, where $(j,j+1)$ are the two branches
merging at the node, \cf. \autoref{fig.beta}.
This fact can easily be seen using a local model. In a
small neighborhood of $x=0$,
one can assume that $\tilde C$ is given by $(y-a)(y-b)=0$. Let
$\tilde\Phi\subset\Sigma_{\tilde d}$ and $\Phi\subset\Sigma_d$
be the disk $\{y=c,\ \ls|x|\le1\}$, $c\in\R$ and $c\gg\ls|a|,\ls|b|$.
Then, the relevant part of~$C$ is the node $(y-ax)(y-bx)=0$,
and $\tilde\Phi$ projects onto the disk
$\tilde\Phi'=\{y=cx,\ \ls|x|\le1\}$, which meets~$\Phi$ at
$(1,c)$. Now, consider one full turn $x=\exp(2\pi it)$, $t\in[0,1]$, and
follow the point $(x,cx)$ in $\partial\tilde\Phi'$: it describes the circle
$y=c\exp(2\pi it)$ encompassing once the two points of
the intersection
%of~$C$ with the fiber.
$C\cap F^\circ_1$.
The class $\Ga_j\Ga_{j+1}$ of this circle
%which is $\Ga_j\Ga_{j+1}$,
is the slope.
Even more precisely, one should start with the constant path $[0,1]\to(1,c)$
%in the fiber
and homotope this path in $F^\circ_x\sminus C$,
keeping one end in~$\Phi$ and the other, in~$\tilde\Phi'$. In the terminal
position, the path is a loop again, and its class
$\Ga_j\Ga_{j+1}$ is the slope.
\endexample

Define the
\emph{twisted local monodromy}
$\tlm_i:=\lm_i\cdot\inn\slope_i$,
%\ie,
where $\inn\:G\to\Aut G$ is the homomorphism sending an element $g$ of a
group~$G$ to the inner automorphism $\inn g\:h\mapsto g\1hg$.
Thus,
$\tlm_i\:\pi(x_i^+)\to\pi(x_i^+)$ is the map
$\Ga\mapsto\slope_i\1(\Ga\ra\lm_i)\slope_i$.
In general, $\tlm_i$ is not a braid.
Take $x_0^+=x_\infty^+$ for the reference fiber and consider the braids
\[*
\pth_i:=\prod_{j=1}^i\Gg_{j-1}\Gb_j\:\pi(x_0^+)\to\pi(x_i^+),\quad
i=1,\ldots,r+1=\infty
\]
(left to right product), the (global) slopes
$\tslope_i:=\slope_i\ra\pth_1\1\in\pi(x_0^+)$, $i=1,\ldots,r$,
and the twisted monodromy homomorphisms
\[*
\tbm_i:=\pth_i\tlm_i\pth_i\1\:\pi(x_0^+)\to\pi(x_0^+),\quad
i=1,\ldots,r.
\]
The following theorem is essentially due to Zariski and van
Kampen~\cite{vanKampen}, and the particular case of improper curves in Hirzebruch
surfaces,
treated by means of the slopes, is considered
in details
in~\cite[\S5.1.3]{degt:book}. Here, we state and outline the proof of
a very special case of this
approach, incorporating the (partial) computation of the braid monodromy of
a real tetragonal curve in terms of its real part.

We use the following common convention:
given an automorphism $\Gb$ of the free group $\<\Ga_1,\ldots,\Ga_4\>$, the
\emph{braid relation $\Gb=\id$}
stands for the quadruple of relations $\Ga_j\ra\Gb=\Ga_j$, $j=1,\ldots,4$.
Note that, since $\Gb$ is an automorphism,
this is equivalent to the infinitely many relations
$\Ga=\Ga\ra\Gb$, $\Ga\in\<\Ga_1,\ldots,\Ga_4\>$.

%\theorem[Zariski--van Kampen]\label{th.vanKampen}
\theorem\label{th.vanKampen}
In the notation above,
the inclusion
of a the reference fiber
induces an epimorphism
$\pi(x_0^+)=\<\Ga_1,\ldots,\Ga_4\>\onto\pi_1$,
and the relations $\tbm_i=\id$, $i=1,\ldots,r$, hold in~$\pi_1$.
If the fiber $x=\infty$ is nonsingular and all non-real singular fibers are
proper for~$\tilde C$, then one also has the
\emph{relations at infinity}
$\pth_{\infty}=\id$ and
$(\Ga_1\ldots\Ga_4)^d=\tslope_r\ldots\tslope_1$.
If, in addition, $C$ has at most one pair of conjugate non-real singular
fibers, then the relations listed define~$\pi_1$.
\endtheorem

\proof
The assertion is a restatement of the classical Zariski--van Kampen theorem
modified for the case of improper curves,
see \cite[Theorem 5.50]{degt:book}.
The relation at infinity
$(\Ga_1\ldots\Ga_4)^d=\tslope_r\ldots\tslope_1$ holds in~$\pi_1$
whenever all slopes not accounted for, namely those at the non-real fibers,
are known to be trivial.
The automorphism
$\pth_{r+1}\:\pi(x_0^+)\to\pi(x_{r+1}^+)=\pi(x_0^+)$ is the monodromy
along the `boundary' of the upper half-plane $\Im x>0$, see
\autoref{fig.base}, \ie, the product of the monodromies about all singular
fibers in this half-plane; if the slopes at these fibers are all trivial,
then $\pth_{r+1}=\id$ in~$\pi_1$.
Finally, if $\tilde C$ has at most one pair of
conjugate non-real singular fibers, then
all but possibly one braid relations are present and hence they define the
group, see \cite[Lemma 5.59]{degt:book}.
\endproof

\section{The computation}\label{S.proofs}

\subsection{The strategy}\label{s.strategy}
We start with a plane sextic $D\subset\Cp2$ and choose homogeneous coordinates
$(z_0:z_1:z_2)$ so that $D$ has a singular point of type~$\bA_\Mn$,
$\Mn\ge3$, at $(0:0:1)$ tangent to the axis $\{z_1=0\}$.
Then, in the affine coordinates $x:=z_1/z_0$, $y:=z_2/z_0$, the curve~$D$ is
given by a polynomial~$\tilde f$ as in \autoref{rem.slopes}, and the same
polynomial~$\tilde f$ defines a
certain tetragonal curve $\tilde C\subset\Sigma_1$,
\viz. the proper transform of~$D$ under the blow-up of $(0:0:1)$.
The common fundamental group
\[*
\pi_1:=\pi_1(\Cp2\sminus D)=\pi_1(\Sigma_1\sminus(\tilde C\cup E))
\]
is computed using \autoref{th.vanKampen} applied to~$\tilde C$ and its
transform $C\subset\Sigma_2$,
with the only nontrivial slope
$\slope=\Ga_1\Ga_2$ or $\Ga_3\Ga_4$ over
$x=0$ given by \autoref{ex.slope}.
(Here, $E\subset\Sigma_1$ is the exceptional section, \iq. the exceptional
divisor over the point $(0:0:1)$ blown up.)
\latin{A priori}, \autoref{th.vanKampen} may only produce a certain
group~$\ttg$ that surjects onto~$\pi_1$ rather than $\pi_1$ itself;
however, in most cases this
group~$\ttg$ is `minimal expected' (\cf. \autoref{s.torus} below) and we do
obtain~$\pi_1$.

The assumption that the fiber $x=\infty$ is nonsingular is not essential as
long as the singularity over~$\infty$ is taken into consideration:
one can always move~$\infty$ to a generic point by a
real projective change of
coordinates. To keep the defining equations as simple as possible, we assume
such a change of coordinates implicitly.
Furthermore, it is only the cyclic order of the singular fibers in the circle
$\Rp1$ that matters, and sometimes we reorder the fibers by applying a cyclic
permutation to their `natural' indices.
In other words, the braid $\Gb_\infty=\Delta^2$ is in the center of~$\BG4$
and, hence, it can be inserted at any place in the relation
$\Gg_0\Gb_1\Gg_1\ldots\Gg_r\Gb_\infty=\id$.

To compute the braids, we outline the real (bold lines) and imaginary
(dotted lines) branches of~$C$ in the figures. Recall that it
is only the mutual position of the real branches and their intersection
indices with the imaginary ones that matters, see \autoref{rem.gamma}.
The `special' node that contributes the only non-trivial slope
(the blow-up center in the passage from~$C$ to~$\tilde C$,
see \autoref{rem.slopes})
is marked with a white dot;
the other singular points of~$C$ (including those of type~$\bA_0$)
are marked with black dots.
The shape of the curve can mostly be recovered using Remarks~\ref{rem.q}
and~\ref{rem.g3}; however, it is usually easier to determine the mutual
position of the roots directly \via\ \Maple.
The braids~$\Gb_i$, $\lm_i$, and $\Gg_j$ are computed from the figures as
explained in \autoref{s.monodromy}.

\warning
The polynomial~$f\nr$ given by \autoref{rem.slopes} is used to determine the
slope and mutual position of the two singular points over $x=0$: the
`special' node is always at $(0,0)$.
For all other applications, \eg, for \autoref{prop.gamma}, this polynomial
should be converted to the reduced form~\eqref{eq.f}.
\endwarning

\subsection{Relations}\label{s.relations}
Recall that a braid relation $\tbm_i=\id$ stands for a quadruple of relations
$\Ga_j\ra\tbm_i=\Ga_j$, $j=1,\ldots,4$. Alternatively, this can be regarded
as an infinite sequence of relations $\Ga\ra\tbm_i=\Ga$, $\Ga\in\FG4$, or,
equivalently, as a quadruple of relations
$\Ga'_j\ra\tbm_i=\Ga'_j$, $j=1,\ldots,4$, where $\Ga_1',\ldots,\Ga_4'$ is
\emph{any} basis for~$\FG4$.
For this reason, in the computation below we start with the braid relations
$\Ga_j'\ra\tlm_i=\Ga_j'$ in the canonical basis over~$x_i^+$ and translate
them to~$x_0^+$ \via~$\pth_i\1$.
In the most common case $\tlm_i=\Gs_r^p$, $r=1,2,3$, $p\in\Z$,
the whole quadruple is equivalent
to the single relation $\{\Ga'_r,\Ga'_{r+1}\}_p=1$, where
%$\{\Ga,\Gb\}_{2k}:=(\Ga\Gb)^k(\Gb\Ga)^{-k}$ and
%$\{\Ga,\Gb\}_{2k+1}:=(\Ga\Gb)^k\Ga(\Ga\Gb)^{-k}\Gb\1$.
\[*
\{\Ga,\Gb\}_{2k}:=(\Ga\Gb)^k(\Gb\Ga)^{-k},\quad
\{\Ga,\Gb\}_{2k+1}:=(\Ga\Gb)^k\Ga(\Ga\Gb)^{-k}\Gb\1.
\]

\remark
\def\Gap{\Ga'}
The braid relations about the fiber $x_k=0$ with the only nontrivial slope,
see \autoref{ex.slope}, can also be presimplified.
Let $\Gap_1,\ldots,\Gap_4$ be the canonical basis in $x_k^+$.
If $\slope_k=\Gap_1\Gap_2$ and $\lm_k=\Gs_1^2\Gs_3^p$,
the braid relations $\tlm_k=\id$ and relation at infinity
$(\Gap_1\ldots\Gap_4)^2=\slope_k$ together are equivalent to
\[*
\Gap_1\Gap_2(\Gap_3\Gap_4)^2=\{\Gap_3,\Gap_4\}_{p+4}=1.
\]
Similarly, if $\slope_k=\Gap_3\Gap_4$ and $\lm_k=\Gs_1^p\Gs_3^2$,
we obtain
\[*
(\Gap_1\Gap_2)^2\Gap_3\Gap_4=\{\Gap_1,\Gap_2\}_{p+4}=1.
\]
Certainly, these relations should be translated back to~$x_0^+$
\via~$\pth_k\1$. Note, though, that we do not use this simplification in the
sequel.
\endremark

\remark\label{rem.xi}
In some cases,
simpler relations are obtained if another point~$x_i^+$, $i>0$, is
taken for the reference fiber. To do so, one merely replaces the
braids~$\pth_j$, $j=1,\ldots,r+1=\infty$, with
$\pth_j':=\pth_i\1\pth_j$.
\endremark

All computations below were performed using \GAP~\cite{GAP4}, with the help
of the
simple braid manipulation routines contained in~\cite{degt:book}.
The \GAP\ code can be downloaded from
\url{http://www.fen.bilkent.edu.tr/~degt/papers/papers.htm}.
The processing is almost fully automated, the input being the braids
$\Gb_i$, $\lm_i$, $\Gg_j$ and the only nontrivial slope $\slope_k=\Ga_1\Ga_2$
or $\Ga_3\Ga_4$, which are read off from the diagrams depicting the curves.

\subsection{The set of singularities \lineref{3A6+A1}}\label{s.3A6+A1}
Any sextic with this set of singularities is $\DG{14}$-special,
see~\cite{degt:Oka}, and,
according to~\cite{degt.Oka}, any $\DG{14}$-special sextic can be given by an
equation of the form
\[*
{\deMaple\obeylines\def
{\\}%
\gathered%
2*t*(t^3-1)*(z[0]^4*z[1]*z[2]+z[1]^4*z[2]*z[0]+z[2]^4*z[0]*z[1])
+(t^3-1)*(z[0]^4*z[1]^2+z[1]^4*z[2]^2+z[2]^4*z[0]^2)#
+t^2*(t^3-1)*(z[0]^4*z[2]^2+z[1]^4*z[0]^2+z[2]^4*z[1]^2)
+2*t*(t^3+1)*(z[0]^3*z[1]^3+z[1]^3*z[2]^3+z[2]^3*z[0]^3)#
+4*t^2*(t^3+2)*(z[0]^3*z[1]^2*z[2]+z[1]^3*z[2]^2*z[0]+z[2]^3*z[0]^2*z[1])
+2*(t^6+4*t^3+1)*(z[0]^3*z[1]*z[2]^2+z[1]^3*z[2]*z[0]^2+z[2]^3*z[0]*z[1]^2)#
+t*(t^6+13*t^3+10)*z[0]^2*z[1]^2*z[2]^2,#
\endgathered}\label{eq.3A6+A1}%
\]
$t^3\ne1$. The set of singularities of this curve is \singset{3A6+A1} if
and only if $t^3=-27$; we use the real value $t=-3$.
After the substitution $z_0=1$, $z_1=x+\frac13$, and $z_2=y/x$ the equation
is brought to the form considered in \autoref{rem.slopes}.
Up to a positive factor, the discriminant~\eqref{eq.discrim}
with respect to~$y$ is
\[*
\deMaple-x^5*(27*x^3-648*x^2+6363*x+7)*(3*x-2)^2*(3*x+1)^7,
\]
which has real roots
\[*
x_1=-\frac13,\quad
x_2\approx-0.001,\quad
x_3=0,\quad
x_4=\frac23,\quad
x_5=\infty
\]
and two simple imaginary roots.
Hence, \autoref{th.vanKampen} does compute the group.

\figure[ht]
\centerline{\cpic{3A6+A1}}
\caption{The set of singularities \lineref{3A6+A1}}\label{fig.3A6+A1}
\endfigure
The only root of~$q$ on the real segment $[-\infty,x_1]$ is $x'\approx-3.48$,
and $b_1(x')<0$; hence, one has $\Gg_0=\id$, see \autoref{prop.gamma}.
The other braids~$\Gb_i$,
$\Gg_j$ are easily found from \autoref{fig.3A6+A1},
and, using \autoref{th.vanKampen} and \GAP, we obtain a group of order~$42$.
This concludes the proof of \autoref{th.3A6+A1}.
\qed

\subsection{Sextics of torus type}\label{s.torus}
All maximal, in the sense of degeneration, sextics of torus type
are described in~\cite{Oka.Pho:moduli}, where a sextic~$D$ is represented by a
pair of polynomials $f_2(x,y)$, $f_3(x,y)$ of degree~$2$ and~$3$,
respectively, so that the defining polynomial of~$D$ is
$f\trs:=f_2^3+f_3^2$.
(Below, these equations are cited in a slightly simplified form: I tried to
clear the denominators by linear changes of variables and appropriate
coefficients.)
Each curve (at least, each of those considered below) has a type~$\bA_\Mn$,
$\Mn\ge3$, singularity at $(0,0)$ tangent to the $y$-axis.
Hence, we start with the substitution
$\tilde f(x,y):=y^6f\trs(x/y,1/y)$ to obtain a polynomial $\tilde f$ as in
\autoref{rem.slopes}; then we proceed as in \autoref{s.strategy}.

To identify the group
$\ttg$ given by \autoref{th.vanKampen}
as~$\MG$, we use the following
\GAP\ code, which was suggested to me by E.~Artal:
\[
\aligned
 &\text{\tt P := PresentationNormalClosure(g, Subgroup(g, a));}\\
 \noalign{\vspace{-3pt}}
 &\text{\tt SimplifyPresentation(P);}
\endaligned
\label{eq.GAP}
\]
here, $\tta$ is an appropriate ratio $\Ga_i\Ga_j\1$
%of two generators
which normally generates the commutant of~$\ttg$.
If the resulting presentation has two generators and no relations,
%then
we conclude that
$\ttg=\pi_1=\MG$,
{\em even when the statement of \autoref{th.vanKampen} does not guarantee a
complete set of relations}.
Indeed, \latin{a priori} we have epimorphisms
$\ttg\onto\pi_1\onto\MG$ (the latter follows from the fact that the curve is
assumed to be of torus type), which induce epimorphisms
$[\ttg,\ttg]\onto[\pi_1,\pi_1]\onto[\MG,\MG]=\FG2$ of the commutants. If
$[\ttg,\ttg]=\FG2$, both these epimorphisms are isomorphisms (since $\FG2$ is
Hopfian) and the $5$-lemma implies that
$\ttg\onto\pi_1\onto\MG$ are also isomorphisms.

In fact, in some cases
(\eg, in \autoref{s.(A14+A2)+A3} and \autoref{s.(A14+A2)+A2+A1}),
the call
{\tt SimplifiedFpGroup(g)}
returns a recognizable presentation of~$\MG$.

\subsection{The set of singularities \lineref{(A14+A2)+A3}}\label{s.(A14+A2)+A3}
The curve in question is \nt{139} in \cite{Oka.Pho:moduli}:
\[*
\deMaple
\aligned
f_2&=80*(-36*y^2+120*x*y-82*x^2+2*x),\\
f_3&=100*(-1512*y^3+7794*y^2*x-18*y^2-11664*y*x^2+144*x*y+5313*x^3-194*x^2+x).
\endaligned
\]
Up to a positive coefficient, the discriminant of~$f\nr$ is
\[*
\deMaple
x^{13}*(5120*x^4+36864*x^3+3456*x^2-2160*x-405)*(x-1)^3.
\]
It has five real roots, which we reorder cyclically as follows:
\[*
x_1=0,\quad
x_2\approx0.27,\quad
x_3=1,\quad
x_4=\infty,\quad
x_5\approx-7.1.
\]
Besides, there are two conjugate imaginary singular fibers, which are of
type~$\bA_0$.

\figure[ht]
\centerline{\cpic{A14+A2+A3}}
\caption{The set of singularities \lineref{(A14+A2)+A3}}\label{fig.(A14+A2)+A3}
\endfigure
The curve is depicted in \autoref{fig.(A14+A2)+A3}, from which all
braids~$\Gb_i$, $\Gg_j$ are easily found.
Taking~$x_0^+$ for the reference fiber and using $\tta=\Ga_1\Ga_2\1$
in~\eqref{eq.GAP}, we obtain $\pi_1=\MG$.

\subsection{The set of singularities \lineref{(A14+A2)+A2+A1}}\label{s.(A14+A2)+A2+A1}
The curve is \nt{142} in \cite{Oka.Pho:moduli}:
\[*
\deMaple
\aligned
f_2&=-45*y^2-240*y*x-106*x^2+90*x,\\
f_3&=1025*y^3+6045*y^2*x-375*y^2+5490*y*x^2-4050*y*x+1354*x^3-2040*x^2+750*x.
\endaligned
\]
Up to a positive coefficient, the discriminant of $f\nr$ is
\[*
\deMaple
x^{13}*(8*x^3-10720*x^2+14250*x-5625)*(x+1)^2*(14*x+15)^3,
\]
and all its roots are real:
\[*
x_1=-\frac{15}{14},\quad
x_2=-1,\quad
x_3=0,\quad
x_4\approx1338,\quad
x_5=\infty.
\]
\figure[ht]
\centerline{\cpic{A14+A2+A2+A1}}
\caption{The set of singularities \lineref{(A14+A2)+A2+A1}}\label{fig.(A14+A2)+A2+A1}
\endfigure
The braids~$\Gb_i$, $\Gg_j$ are found from \autoref{fig.(A14+A2)+A2+A1} and,
using~$x_0^+$ as the reference fiber and $\tta=\Ga_1\Ga_2\1$
in~\eqref{eq.GAP}, we conclude that $\pi_1=\MG$.

\subsection{The set of singularities \lineref{(A11+2A2)+A4}}\label{s.(A11+2A2)+A4}
This is \nt{118} in \cite{Oka.Pho:moduli}:
\[*
\deMaple
\aligned
f_2&=\frac15(-3456*y^2+1200*y*x-3005*x^2+240*x),\\
f_3&=\frac15(-89856*y^3+130464*y^2*x-6912*y^2-112680*y*x^2+8640*y*x+\\
\noalign{\vspace{-8pt}}
&\phantom{{}=\dfrac15(}91345*x^3-13320*x^2+480*x).
\endaligned
\]
Up to a positive coefficient, the discriminant of~$f\nr$ is
\[*
\deMaple
-x^{10}*(25*x^3+290*x^2+360*x+162)*(35*x^2-384*x+1152)^3.
\]
It has three real roots, which we reorder cyclically as follows:
\[*
x_1=0,\quad
x_2=\infty,\quad
x_3\approx-10.26.
\]
\figure[ht]
\centerline{\cpic{A11+2A2+A4}}
\caption{The set of singularities \lineref{(A11+2A2)+A4}}\label{fig.(A11+2A2)+A4}
\endfigure
In addition, there are two pairs of complex conjugate singular fibers, of
types~$\bA_2$ and~$\bA_0$. Thus, \latin{a priori} \autoref{th.vanKampen} only
gives us a certain epimorphism $\ttg\onto\pi_1$.
However, using $\tta=\Ga_1\Ga_2\1$ in~\eqref{eq.GAP}, we conclude that
$\ttg=\pi_1=\MG$. (All braids are found from \autoref{fig.(A11+2A2)+A4} and
the reference fiber is~$x_1^+$, see \autoref{rem.xi}.)

\subsection{The set of singularities \lineref{(A8+3A2)+A4+A1}}\label{s.(A8+3A2)+A4+A1}
This curve is \nt{83} in \cite{Oka.Pho:moduli}:
\[
{\deMaple
\aligned
f_2&=-565*y^2-14*y*x+176*y-5*x^2+104*x-16,\\
f_3&=13321*y^3+3135*y^2*x-6294*y^2+207*y*x^2-3516*y*x+1056*y+\\
\noalign{\vspace{-2pt}}
&\hphantom{{}={}}25*x^3-558*x^2+624*x-64.
\endaligned}
\label{eq.(A8+3A2)+A4+A1}
\]
Up to a positive coefficient, the discriminant of~$f\nr$ is
\[*
\deMaple
x^3*(x+3)*(x+9)^2*(11915*x^3+96579*x^2-14823*x+729)^3*(x-9)^9.
\]
It has five real roots, which we reorder cyclically as follows:
\[*
x_1=0,\quad
x_2=9,\quad
x_3=-9,\quad
x_4\approx-8.26,\quad
x_5=-3.
\]
We conclude that the curve has only two non-real singular fibers, which are
cusps.
Hence, \autoref{th.vanKampen} gives us a complete presentation of~$\pi_1$.

\figure[ht]
\centerline{\cpic{A8+3A2+A4+A1}}
\caption{The set of singularities \lineref{(A8+3A2)+A4+A1},
projected from~$\bA_4$}\label{fig.(A8+3A2)+A4+A1}
\endfigure
In the interval $(x_5,x_1)$, where $f$ has four imaginary branches, $q$ has
four roots
\[*
x_1'\approx-2.93,\quad
x_2'=-1.92,\quad
x_3'\approx-0.79,\quad
x_4'\approx-0.14,
\]
with $b_1$ negative at~$x_1'$,~$x_3'$ and positive at~$x_2'$,~$x_4'$;
at the latter two points one also has $q'<0$. Hence, $\Gg_0=\tau^{-2}$, see
\autoref{prop.gamma}.
All other braids are esily found from \autoref{fig.(A8+3A2)+A4+A1}.

\remark\label{rem.rho.inf}
For a further simplification, observe that
the braid $\pth_\infty$ appearing in \autoref{th.vanKampen} equals
\[*
\Gs_2\1\Gs_1\Gs_3\1\Gs_1\Gs_3\1\Gs_2\cdot
\Gs_1\1\cdot
\Gs_2\1\Gs_1\cdot
\Gs_2^{-4}\cdot
\Gs_3\1\cdot
\Gs_2^{-2}\cdot
\Gs_3\1\Gs_2\cdot
\Gs_1\1\cdot
(\Gs_3\Gs_1\Gs_2)^4,
\]
and one can check that $\pth_\infty=\pim\1\Gs_1^3\pim$, where
$\pim:=\Gs_2\Gs_1\1\Gs_3^2\Gs_2$.
(Note that $\pth_\infty$ represents the monodromy about a single imaginary
cusp of the curve; hence, it
%\emph{must}
is expected to
be conjugate to~$\Gs_1^3$.)
Thus, we can replace the quadruple of relations $\pth_\infty=\id$ with a
single relation $\{\Ga_1,\Ga_2\}_3\ra\pim=1$, \cf. \autoref{s.relations}.
\endremark

Now, taking $x_3^+$ for the reference fiber, see \autoref{rem.xi},
using \autoref{rem.rho.inf},
and applying
%\GAP\ function
{\tt SimplifiedFpGroup(g)}, we arrive at~\eqref{gr.(A8+3A2)+A4+A1}. This
presentation has
three generators and
four relations of total length $48$.
Together with the previous sections, this concludes the proof of
\autoref{th.torus}.
\qed

\remark\label{rem.unknown}
The Alexander module of the group~$\pi_1$ considered in this section is
$\Z[t,t\1]/(t^2-t+1)$,
and the finite quotients $\pi_1/\Ga_2^p$, $p=2,3,4$, are
isomorphic to the similar quotients of~$\MG$.
My laptop failed to compute the order of $\pi_1/\Ga_2^5$.
\endremark

\remark\label{rem.A8}
In \eqref{eq.(A8+3A2)+A4+A1}, the singular point at the origin is of
type~$\bA_4$. One can start with a change of variables
\figure[ht]
\centerline{\cpic{A8+3A2+A4+A1-alt}}
\caption{The set of singularities \lineref{(A8+3A2)+A4+A1},
projected from~$\bA_8$}\label{fig.(A8+3A2)+A4+A1-alt}
\endfigure
$x\mapsto y+9$, $y\mapsto x+1$ and resolve the type~$\bA_8$ point instead.
The
%curve
tetragonal model
is depicted in \autoref{fig.(A8+3A2)+A4+A1-alt}, and
the computation becomes slightly simpler, but the resulting presentation is
of the same complexity, even
%if one observes that
with the additional observation that
$\pth_\infty=\pim\1\Gs_1^3\pim$, where
$\pim:=\Gs_2\Gs_1\1\Gs_3\Gs_2$, \cf. \autoref{rem.rho.inf}.
\endremark

\subsection{Proof of \autoref{prop.pert}}\label{proof.pert}
For the sets of singularities
\lineref{(A14+A2)+A3},
\lineref{(A14+A2)+A2+A1}, and
\lineref{(A11+2A2)+A4},
the statement
%follows from
is an immediate consequence of
\cite[Theorem 7.48]{degt:book}.
For \lineref{3A6+A1}, the only proper quotient of the commutant
$[\pi_1,\pi_1]=\CG7$ is trivial; hence, the group~$\pi_1'$ of any
perturbation~$D'$ is either abelian, $\pi_1'=\CG6$, or isomorphic to~$\pi_1$,
the latter being the case if and only if $D'$ is $\DG{14}$-special,
see~\cite{degt:Oka}.

For the remaining set of singularities \lineref{(A8+3A2)+A4+A1},
we proceed as follows. Any proper perturbation
factors through a maximal one, where a single singular point~$P$ of
type~$\bA_\Mn$ splits into two points~$\bA_{\Mn'}$, $\bA_{\Mn''}$, so that
$\Mn'+\Mn''=\Mn-1$.
Assume that $P\ne(0:0:1)$, see \autoref{s.strategy}.
Then
this point
corresponds to a certain singular fiber~$x_i$ of
the tetragonal model~$C$ and
gives rise to a braid relation $\{\Ga_k,\Ga_{k+1}\}_{\Mn+1}\ra\pth_i\1=1$, see
\autoref{s.relations}.
For the new curve~$D'$, this relation changes to
$\{\Ga_k,\Ga_{k+1}\}_s\ra\pth_i\1=1$, where $s:=\gcd(\Mn'+1,\Mn''+1)$.

For any perturbation of any point~$P$, we have $s=3$
if $P$ is of type~$\bA_8$ or~$\bA_2$ and the result is still of torus type,
and $s=1$ otherwise.
Now, the statement is easily proved by repeating the computation with the
braid $\lm_i=\Gs_k^{\Mn+1}$ replaced with $\Gs_k^s$.
(If it is the type~$\bA_4$ point that is perturbed,
one can use the alternative tetragonal model given by
\autoref{rem.A8}.)
\qed

{
\let\.\DOTaccent
\def\cprime{$'$}
\bibliographystyle{amsplain}
\bibliography{degt}
}

\end{document}